\input amstex
\magnification=\magstep1
\input epsf
\baselineskip=13pt
\documentstyle{amsppt}
\vsize=8.7truein
\CenteredTagsOnSplits
\NoRunningHeads
\def\today{\ifcase\month\or
  January\or February\or March\or April\or May\or June\or
  July\or August\or September\or October\or November\or December\fi
  \space\number\day, \number\year}

\def\conv{\operatorname{conv}}

\def\Sym{\operatorname{Sym}}

\topmatter
\title Approximating a Norm by a Polynomial \endtitle
\author Alexander Barvinok \endauthor
\address Department of Mathematics, University of Michigan, Ann Arbor,
MI 48109-1109 \endaddress
\email barvinok$\@$umich.edu \endemail
\date  \today \enddate
\thanks This research was partially supported by NSF Grant DMS 9734138.
\endthanks
\abstract We prove that for any norm $\| \cdot \|$ in the 
$d$-dimensional real vector space $V$ and for any odd $n>0$ there 
is a non-negative polynomial $p(x)$, $x \in V$ of degree $2n$ such that 
$$p^{1\over 2n}(x) \leq \|x\| \leq {n+d-1 \choose n}^{1 \over 2n}
 p^{1\over 2n}(x).$$
Corollaries and polynomial approximations of the Minkowski functional 
of a convex body are discussed.   
\endabstract
\keywords norm, approximation, polynomial, John's ellipsoid, computational
complexity
\endkeywords 
\endtopmatter
\document

\head 1. Introduction and the Main Result \endhead

Our main motivation is the following general question.
Let us fix a norm $\| \cdot \|$ in a finite dimensional real 
vector space $V$ (or, more generally, the Minkowski functional of a 
convex body in $V$). Given a point $x \in V$, how fast can one 
compute or approximate $\|x\|$? For example, various optimization problems
can be posed this way.
 As is well known, (see, for example, 
Lecture 3 of [1]), any  norm in $V$ can be approximated by an $\ell^2$
norm in $V$ within a factor of $\sqrt{\dim V}$. From the computational 
complexity point of view, an $\ell^2$ norm of $x$ is just the square root
of a positive definite quadratic form $p$ in $x$ and hence can be computed 
``quickly'', that is, in time polynomial in $\dim V$ for 
any $x \in V$ given by its coordinates in some basis of $V$. Note, that 
we do not count the time required for ``preprocessing'' the norm 
to obtain the quadratic form $p$, as we consider the norm fixed 
and not a part of the input.  
It turns out that by employing higher degree forms $p$, we can  
improve the approximation: for any $c>0$, given an 
$x \in V$, one can approximate $\|x\|$ within a factor of $c\sqrt{\dim V}$ 
in time polynomial in $\dim V$. This, and some other approximation results 
follow easily from our main theorem. 
\proclaim{(1.1) Theorem} Let $V$ be a $d$-dimensional real vector space 
and let $\| \cdot \|:$ $V \longrightarrow {\Bbb R}$ be a norm in $V$. 
For any odd integer $n >0$ there exists a 
homogeneous polynomial $p: V \longrightarrow 
{\Bbb R}$ of degree $2n$ such that $p(x) \geq 0$ and 
$$p^{1 \over 2n}(x) \leq \| x\| \leq {n+d-1 \choose n}^{1 \over 2n}
 p^{1 \over 2n}(x)$$
for all $x \in V$.
\endproclaim
We prove Theorem 1.1 in Section 2.

Let us fix an $n$ in Theorem 1.1. Then, as $d$ grows, the value 
of $\displaystyle p^{1 \over 2n}(x)$ approximates $\| \cdot \|$ within 
a factor of $c_n \sqrt{d}$, where 
$\displaystyle c_n \approx (n!)^{-{1 \over 2n}} \approx \sqrt{e/n}$.
Since for any fixed $n$, computation of $p(x)$ takes a $d^{O(n)}$ time,
for any $c>0$ we obtain a polynomial time algorithm to approximate 
$\|x\|$ within a factor of $c \sqrt{d}$ (again, we do not count the time 
required for preprocessing, that is, to find the polynomial $p$).

If we allow $n$ to grow linearly with $d$, we can get a constant factor 
approximation. Indeed, if we choose $n=\gamma d$ for some 
$\gamma>0$, for large $d$ we have 
$${n+d-1 \choose n}^{1 \over 2n} \approx 
\exp\Bigl\{{1 \over 2} \ln {\gamma+1\over\gamma}+ 
{1 \over 2\gamma} \ln (\gamma+1)
\Bigr\}.$$
Since for any fixed $\gamma>0$, computation of $p(x)$ takes $2^{O(d)}$ time,
for any $c>1$ we can get an algorithm of $2^{O(d)}$ complexity approximating 
the value of $\|x\|$ within a factor of $c$.

One can hope that for special norms $\| \cdot \|$ (for example, ones 
with a large symmetry group) one can obtain better 
approximability/computability 
results due to special features of the polynomials $p$ (for example, 
invariance with respect to the action of a large symmetry group).
Indeed, the construction of the proof of Theorem 1.1 (see 
Section 2) preserves, for example, group invariance.

\head 2. Proof of Theorem 1.1 \endhead

Let $B$ be the unit ball of $\| \cdot \|$, so 
$$B=\bigl\{x \in V: \|x\| \leq 1 \bigr\}.$$
Hence $B$ is a centrally symmetric convex compact set containing 
the origin in its interior.

Let $V^{\ast}$ be the dual space of all 
linear functions $f: V \longrightarrow {\Bbb R}$ and let $C \subset V^{\ast}$ 
be the polar of $B$:
$$C=\Bigl\{f \in V^{\ast}: f(x) \leq 1 \quad 
\text{for all} \quad x \in B \Bigr\}.$$
Hence $C$ is a centrally symmetric convex compact set containing the 
origin in its interior.
Using the standard duality argument (see, for example, Section 1.6 of
[3]), we can write
$$\|x\|=\max_{f \in C} f(x). \tag2.1$$
Let 
$$W=V^{\otimes n}=\underbrace{V \otimes \ldots \otimes V}_{\text{$n$ times}}
\quad \text{and} \quad W^{\ast}=\bigl(V^{\otimes n}\bigr)^{\ast}=
\underbrace{V^{\ast} \otimes \ldots \otimes V^{\ast}}_{\text{$n$ times}}$$
be the $n$-th tensor powers of $V$ and $V^{\ast}$ respectively.

For vectors $x \in V$ and $f \in V^{\ast}$ let 
$$x^{\otimes n}=\underbrace{x \otimes \ldots \otimes x}_{\text{$n$ times}}
\quad \text{and} \quad 
f^{\otimes n}=\underbrace{f \otimes \ldots \otimes f}_{\text{$n$ times}}$$
denote the $n$-th tensor power $x^{\otimes n} \in W$ and 
$f^{\otimes n} \in W^{\ast}$ respectively.
 
By (2.1), we can write 
$$\|x\|^n=\max_{f \in C} \bigl(f(x)\bigr)^n=
\max_{f \in C} f^{\otimes n}(x^{\otimes n}).
\tag2.2$$
Let $D$ be the convex hull of $f^{\otimes n}$ for $f \in C$:
$$D=\conv\bigl\{f^{\otimes n}: f \in C \bigr\}.$$
Then $D$ is a convex compact centrally symmetric (we use that $n$ is odd) 
subset of $W^{\ast}$ and from (2.2) we can write 
$$\|x\|^n=\max_{f \in C} f^{\otimes n}(x^{\otimes n})=
\max_{g \in D} g(x^{\otimes n}). \tag2.3$$ 
Let us estimate the dimension of $D$. There is a natural action of the 
symmetric group $S_n$ in $W^{\ast}$ which permutes the factors 
$V^{\ast}$, so that 
$$\sigma(f_1 \otimes \ldots \otimes f_n)=f_{\sigma^{-1}(1)} \otimes 
\ldots \otimes f_{\sigma^{-1}(n)}.$$
Let $\Sym(W^{\ast}) \subset W^{\ast}$ 
be the {\it symmetric part} of $W^{\ast}$,
that is, the invariant subspace of that action.
As is known, the dimension of $\Sym(W^{\ast})$ is that of the space 
of homogeneous polynomials of degree $n$ in $d$ real variables 
(see, for example,
Lecture 6 of [2]).
Next, we observe that $f^{\otimes n} \in \Sym(W^{\ast})$ for all 
$f \in V^{\ast}$ and, therefore,
$$\dim D \leq \dim \Sym(W^{\ast}) ={n+d-1 \choose n} \tag2.4$$ 
Let $E$ be the John's ellipsoid of $D$ in the affine hull 
of $D$, that is the (unique) ellipsoid 
of the maximum volume inscribed in $D$. As is known, (see, for 
example, Lecture 3 of [1]) 
$$E \subset D \subset \bigl(\sqrt{\dim D}\bigr) E.$$
Combining this with (2.3), we write 
$$\max_{g \in E} g(x^{\otimes n})\leq 
\|x\|^n \leq \bigl(\sqrt{\dim D} \bigr)\max_{g \in E} g(x^{\otimes n})$$
and, by (2.4),
$$\max_{g \in E} g(x^{\otimes n}) \leq \|x\|^n \leq 
{n+d-1 \choose n}^{1 \over 2} \max_{g \in E} g(x^{\otimes n}).
\tag2.5$$
Let 
$$q(x)=\max_{g \in E} g(x^{\otimes n}).$$
We claim that $p(x)=q^2(x)$ is a polynomial in $x$ of degree $2n$.
Indeed, let us choose a basis $e_1, \ldots, e_d$  
in $V$ and the dual basis $f_1, \ldots, f_d$ in $V^{\ast}$,
so that $f_i(e_j)=\delta_{ij}$. Then 
$W$ acquires the basis 
$$e_{i_1 \ldots i_n}=e_{i_1} \otimes \ldots \otimes e_{i_n} \quad 
\text{for} \quad 1 \leq i_1, \ldots, i_n \leq d$$
and $W^{\ast}$ acquires the basis
$$f_{i_1 \ldots i_n} =f_{i_1} \otimes \ldots \otimes f_{i_n} \quad 
\text{for} \quad 1 \leq i_1, \ldots, i_n \leq d.$$
Geometrically, $V$ and $V^{\ast}$ are identified with ${\Bbb R}^d$ and 
$W$ and $W^{\ast}$ are identified with ${\Bbb R}^{dn}$.
Let $K \subset W^{\ast}$ be the Euclidean unit ball defined by the inequality
$$K=\Bigl\{h \in W^{\ast}: \sum_{1 \leq i_1, \ldots, i_n \leq d} 
h_{i_1 \ldots i_n}^2 \leq 1 \Bigr\},$$
where $h_{i_1 \ldots i_n}$ is the corresponding coordinate of $h$ 
with respect to the basis $\{f_{i_1 \ldots i_n} \}$.
Since $E$ is an ellipsoid, there is a linear transformation $T: W^{\ast} 
\longrightarrow W^{\ast}$ such that $T(K)=E$. 
Let $T^{\ast}: W \longrightarrow W$ be the conjugate linear transformation
and let $y=T^{\ast}\bigl(x^{\otimes n}\bigr)$. Hence the coordinates
$y_{i_1 \ldots i_n}$ of $y$ with respect to the basis 
$\bigl\{e_{i_1 \ldots i_n}\bigr\}$
are polynomials in $x$ of degree $n$.
Then
$$\split q(x) &=\max_{g \in E} g(x^{\otimes n})=
\max_{h \in K} T(h)\bigl(x^{\otimes n}\bigr)=
\max_{h \in K} h\bigl(T^{\ast}(x^{\otimes n})\bigr) \\ &=\max_{h \in K}
h(y)=
\sqrt{\sum_{1 \leq i_1, \ldots, i_n \leq d} 
y_{i_1 \ldots i_n}^2}. \endsplit$$ 
Hence we conclude that $p(x)=q^2(x)$ is a homogeneous 
polynomial in $x$ of degree $2n$, which is non-negative for 
all $x \in V$ (moreover, $p(x)$ is seen to be a sum of squares). 
From (2.5), we conclude that 
$$p^{1 \over 2n}(x) \leq \|x\| \leq {n + d-1 \choose n}^{1 \over 2n} 
p^{1 \over 2n}(x),$$
as claimed.

\head 3. An Extension to Minkowski Functionals \endhead

There is a version of Theorem 1.1 for Minkowski functionals of not 
necessarily centrally symmetric convex bodies.
\proclaim{(3.1) Theorem} Let $V$ be a $d$-dimensional real vector space,
let $B \subset V$ be a convex compact set containing the origin in its 
interior and let $\|x\|=\inf \bigl\{\lambda>0: x \in \lambda B \bigr\}$ be 
its Minkowski functional. For any odd integer $n>0$ there exist 
a homogeneous polynomial $p: V \longrightarrow {\Bbb R}$ of degree $2n$ 
and a homogeneous polynomial
$r: V \longrightarrow {\Bbb R}$ of degree $n$ such that $p(x) \geq 0$
and
$$\biggl(r(x) + \sqrt{p(x)} \biggr)^{1 \over n} 
\leq \|x \| \leq \biggl(r(x) + {n+d -1 \choose n} 
\sqrt{p(x)}\biggr)^{1 \over n}$$
for all $x \in V$.
\endproclaim 
\demo{Proof} The proof follows the proof of Theorem 1.1 with 
some modifications. Up to (2.4) no essential changes are needed 
(note, however, that now we have to use that $n$ is odd in (2.2)). Then,
since the set $D$ is not necessarily centrally symmetric, 
 we can only find an ellipsoid $E$ (centered at the origin) of 
$W^{\ast}$ and a point $w \in D$, such that
$$E \subset D-w \subset (\dim D)E,$$
see, for example, Lecture 3 of [1].
Then (2.5) transforms into 
$$\max_{g \in E} g(x^{\otimes n}) \leq \|x\|^n -w(x^{\otimes n}) 
\leq {n+d-1 \choose n} \max_{g \in E} g(x^{\otimes n}).$$
Denoting
$$p(x)=\Bigl(\max_{g \in E} g(x^{\otimes n})\Bigr)^2 \quad 
\text{and} \quad r(x)=w(x^{\otimes n})$$
we proceed as in the proof of Theorem 1.1.
{\hfill \hfill \hfill} \qed
\enddemo

\head References \endhead 

\item{1.} K. Ball, An elementary introduction to modern convex geometry,
in: {\it Flavors of Geometry}, 1--58, Math. Sci. Res. Inst. Publ., {\bf 31}, 
Cambridge Univ. Press, Cambridge, 1997.

\item{2.} W. Fulton and J.Harris, {\it Representation Theory. A First Course},
Graduate Texts in Mathematics, {\bf 129},
Readings in Mathematics, Springer-Verlag, New York, 1991.

\item{3.} R. Schneider, {\it Convex Bodies: the Brunn-Minkowski Theory},
Encyclopedia of Mathematics and its Applications, {\bf 44},
Cambridge University Press, Cambridge, 1993.

\enddocument
\end